\documentclass[11pt]{article}

\usepackage[a4paper,margin=31mm]{geometry}
\usepackage{amsmath,amssymb,amsthm,mathtools}
\usepackage[T1]{fontenc}
\usepackage{microtype}
\usepackage[colorlinks=true,linkcolor=blue,citecolor=blue,urlcolor=blue]{hyperref}

\newtheorem{conjecture}{Conjecture}[section]

\newtheorem{remark}[conjecture]{Remark}

\newcommand{\Area}{\operatorname{Area}}
\newcommand{\length}{\operatorname{length}}
\newcommand{\Det}{\operatorname{det}}
\newcommand{\cJ}{\mathcal J}
\newcommand{\cN}{\mathcal N}
\newcommand{\RR}{\mathbb R}

\title{\bfseries Corner contributions to Neumann jump determinants:\\
three model calculations and a BFK conjecture}
\author{Victor Kalvin}
\date{}

\begin{document}

\maketitle

\begin{abstract}
We formulate a local corner-factor conjecture for the determinant of the
Neumann jump operator on a piecewise real-analytic cutting curve.  For the mirror
double of a simply connected geodesic polygon with interior angles
$\pi\alpha_1,\ldots,\pi\alpha_N$, the conjectural determinant is
\[
 \Det_{\angle}'\cN
 =
 \frac{\length(\partial P)}2
 \prod_{j=1}^N\alpha_j^{-1/2}.
\]
Here $\Det_{\angle}'$ denotes an intrinsic determinant, still to be
constructed, that is required to satisfy a
Burghelea--Friedlander--Kappeler gluing formula; only its quotient by
$\length(\partial P)$ is purely angle-dependent.
Three model calculations support the conjecture: a flat polygon and
its double give $\frac12\prod_j\alpha_j^{-1/2}$; a spherical spindle
split into congruent lunes gives $1/(2\alpha)$ for two angles
$\pi\alpha$; and every spherical Coxeter triangle with angles
$(\pi/p,\pi/q,\pi/r)$ gives $\frac12\sqrt{pqr}$, including $\sqrt2$
for the octant.  The conjecture also reduces the Dirichlet
determinant of a constant-curvature polygon to that of its closed
double.  Although the singular anomaly formula applies to the double
in general, an explicit evaluation in nonzero curvature requires
solving its uniformization problem.  We discuss connections with the
Neumann jump determinants in the work of Wiegmann--Zabrodin and Wang
and with the recent Grunsky-operator approach to Coulomb gases on
domains with corners.
\end{abstract}

\section{Introduction}
This note provides a calibration for the missing corner extension of
the BFK formula.  The usual formula involves the determinant of the
Neumann jump operator along a smooth cutting curve \cite{BFK}; its
extension to conical surfaces requires the cut to avoid the singular
set \cite{KalvinJFA}.  A cut through conical points has corners and lies
outside this theory.  This obstruction and the need for a corresponding
BFK formula were recorded in the 2023 preprint of
\cite{KalvinASNS}, where the determinant calculations therefore relied
instead on the singular anomaly formulas of \cite{KalvinJFA}.  The
present note identifies an explicit value that any such extension must
reproduce.

All three calculations concern mirror seams: the two pieces are
isometric and interchanged by reflection.  They therefore isolate the
symmetric BFK contribution rather than the full conformal-welding
problem.

The computations suggest that in the
mirror-double setting, after the correct corner
renormalization, the BFK quotient contains the universal local factor
\[
        \frac12\prod_j \alpha_j^{-1/2},
\]
where \(\pi\alpha_j\) are the angles of the polygonal cut.

This factor also provides a benchmark for comparing the spectral BFK
regularization with the Loewner and Grunsky regularizations for curves
with corners, as discussed in Section~7.

Throughout the note, $\Det'\Delta_M$ denotes the zeta-regularized
determinant of the Laplacian on a closed surface $M$, with the zero
eigenvalue omitted.  If $P$ is a surface with boundary, then
$\Det\Delta_{P,D}$ denotes the zeta-regularized determinant of its
Dirichlet Laplacian.  Whenever the metric has conical singularities, the
Friedrichs extensions of the corresponding Laplacians are understood.

Let a closed surface $M$ be cut along a smooth separating curve $\Gamma$
into two surfaces $M_+$ and $M_-$.  The BFK surgery formula
\cite{BFK} can be written as
\begin{equation}\label{eq:smooth-bfk}
 \frac{\Det'\Delta_M}{\Area(M)}
 =
 \Det\Delta_{M_+,D}\,\Det\Delta_{M_-,D}\,
 \frac{\Det'\cN_\Gamma}{\length(\Gamma)}.
\end{equation}
The same BFK identity was proved for Friedrichs Laplacians on surfaces
with conical singularities in
\cite[Proposition~2.8]{KalvinJFA}, under the essential restriction that
the smooth cutting curve avoid the conical points.  The present problem
is precisely outside that setting: on the mirror double of a polygon
the cutting curve passes through the conical points and is itself
nonsmooth there.

The jump operator is
\[
 \cN_\Gamma f=\partial_{n_+}u_++\partial_{n_-}u_-,
\]
where $u_\pm$ are the harmonic extensions of $f$ and $n_\pm$ are the
outward unit normals.  Its kernel consists of the constants.

Formula \eqref{eq:smooth-bfk} is an extremely effective way of passing
between determinants on a closed surface and Dirichlet determinants on
its pieces.  It is much less clear what the formula should mean when the
cutting curve has corners.  The jump operator is then no longer a
classical pseudodifferential operator on a smooth closed curve.  Corner
Mellin symbols enter its local structure, logarithmic terms may enter
heat and resolvent expansions, and the ordinary zeta-function
definition need not apply without modification. A direct approximation of the corners by smooth arcs
does not solve this problem: zeta
determinants are not, in general, continuous under such a singular
approximation.

Several smooth results motivate this question.  Wiegmann and Zabrodin
\cite{WiegmannZabrodin} related a term in the free energy of the Dyson
gas to the spectral determinant of the Neumann jump operator.  Wang
\cite{Wang} expressed the Loewner energy of a smooth loop through a BFK
determinant quotient.  In the multichordal setting, Peltola and Wang
\cite{PeltolaWang} expressed the Loewner potential as a quotient of
Dirichlet Laplacian determinants.  None of these results extends the
spectral jump determinant or the BFK formula to polygonal or conical
cuts; the singular comparison formulas of \cite{KalvinJFA}
also leave boundary corners open.

For Jordan domains with corners, Johansson and Viklund
\cite{JohanssonViklund} found a logarithmic divergence of the truncated
Grunsky determinant with coefficient
\[
 \sum_j\bigl(\alpha_j+\alpha_j^{-1}-2\bigr)
\]
and conjectured that its reduced finite part is related to Laplacian
determinants; their proof uses analyticity of the boundary arcs.
Courteaut, Johansson and Viklund
\cite{CourteautJohanssonViklund} introduced an arc-Grunsky Fredholm
determinant.  These Fredholm determinants are distinct from the
spectral jump determinant considered here: the two sides of a mirror
double are isometric, whereas the planar interior--exterior problem
contains nontrivial welding data.

The question of a polygonal analogue of the Wiegmann--Zabrodin formula
was posed to the author by Paul Wiegmann.  The author discussed it with
Wiegmann in 2023 and again more recently in connection with
\cite{JohanssonViklund,CourteautJohanssonViklund}; through Wiegmann,
the question was also discussed with colleagues in Sweden.  The
present note was also prompted by a recent question of Werner M\"uller
concerning extremal properties of determinants on hyperbolic polygons.

\section{The BFK quotient and the corner-factor conjecture}

Let $P$ be a compact simply connected surface of constant curvature with
piecewise geodesic boundary.  Suppose that its vertices have interior
angles
\[
 \pi\alpha_1,\ldots,\pi\alpha_N,\qquad 0<\alpha_j<2.
\]
Let $\widehat P=P\cup_{\partial P}P$ be the mirror double.  It is a
closed surface of spherical topology with conical angles $2\pi\alpha_j$ at the vertices.  We
use the Friedrichs Laplacian on $\widehat P$ and the Friedrichs
Dirichlet Laplacian on $P$.

Independently of whether a spectral determinant of the Neumann jump
operator on the piecewise real-analytic curve $\partial P$ has been defined,
the BFK quotient
\begin{equation}\label{eq:def-J}
 \cJ_{\rm BFK}(P):=
 \frac{\Det'\Delta_{\widehat P}/\Area(\widehat P)}
      {\bigl(\Det\Delta_{P,D}\bigr)^2}
\end{equation}
is well defined whenever the Laplacian determinants on the right-hand
side are well defined.  For a smooth mirror double, the jump
operator equals twice the Dirichlet-to-Neumann operator of one copy and
\eqref{eq:smooth-bfk} gives
\[
 \cJ_{\rm BFK}(P)=
 \frac{\Det'\cN_{\partial P}}{\length(\partial P)}.
\]
For a simply connected smooth double, the last quotient equals $1/2$, see
\cite[Lemma~2.10 and Section~3.4]{KalvinJFA}; see also
\cite[Lemma~2.2]{KalvinCCM}.  The corresponding determinant of the
Dirichlet-to-Neumann operator on a smooth planar domain was computed in
\cite{EdwardWu}.
\begin{conjecture}[local corner factor]\label{conj:corner}
There is a canonical determinant $\Det_{\angle}'\cN_{\partial P}$ of the
Neumann jump operator on a piecewise real-analytic curve with finitely
many corners for which the BFK formula remains valid.  In the
mirror-double setting above, it satisfies
\begin{equation}\label{eq:corner-conj}
 \frac{\Det_{\angle}'\cN_{\partial P}}{\length(\partial P)}
 =\cJ_{\rm BFK}(P)
 =\frac12\prod_{j=1}^{N}\alpha_j^{-1/2}.
\end{equation}
\end{conjecture}

Only the normalized determinant in \eqref{eq:corner-conj} is local in
the corner angles; the determinant itself retains the length of the
cut.  The notation $\Det_{\angle}'$ is not meant to denote the ordinary
zeta determinant: its intrinsic definition and its identification with
the Laplacian quotient \eqref{eq:def-J} are part of the conjecture.

\begin{remark}
Setting every $\alpha_j=1$ in \eqref{eq:corner-conj} recovers the smooth
value $1/2$.  In spherical topology, the formula also predicts no
residual dependence of the normalized mirror jump determinant on global
moduli.
\end{remark}

\section{A flat polygon and its double}

Let $\widehat P$ be the flat conical sphere obtained by gluing two
copies of a Euclidean polygon $P$ along their boundaries.  Write
\[
 S=\Area(\widehat P)=2A,\qquad A=\Area(P),
\]
and write the interior angles of $P$ as $\pi\alpha_j$.  Thus
\[
 0<\alpha_j<2,\qquad
 |\boldsymbol\alpha|
 :=\sum_{j=1}^N\alpha_j=N-2,
\]
and the cone angles of $\widehat P$ are $2\pi\alpha_j$.

Choose a conformal coordinate $z$ on $\widehat P=\mathbb{CP}^1$ in
which the conical points $p_1,\ldots,p_N$ lie on $\mathbb R$.  They are
the Schwarz--Christoffel prevertices of $P$ and, modulo real M\"obius
transformations, describe its moduli.  In the notation used in
\cite{KalvinCV}, the metric is
\begin{equation}\label{eq:flat-polygon-metric}
 m_{\boldsymbol\alpha,\boldsymbol p}^{S}
 =
 S C_{\boldsymbol\alpha,\boldsymbol p}^{\,2}
 \prod_{j=1}^N|z-p_j|^{2\alpha_j-2}|dz|^2,
\end{equation}
where
\begin{equation}\label{eq:C-alpha-p-flat}
 C_{\boldsymbol\alpha,\boldsymbol p}^{-2}
 =
 \int_{\mathbb C}
 \prod_{j=1}^N|z-p_j|^{2\alpha_j-2}
 \frac{dz\wedge d\bar z}{-2i}.
\end{equation}
The upper and lower half-planes are the two copies of $P$.

Let $\Delta_{\widehat P}$ be the corresponding Friedrichs Laplacian.
Its spectral zeta function satisfies
\begin{equation}\label{eq:zeta-alpha-polygon}
 \zeta_{\boldsymbol\alpha}(0)
 =
 \frac1{12}\sum_{j=1}^N
 \left(\frac1{\alpha_j}-\alpha_j\right)-1.
\end{equation}
We use the angle-parametrized local function from
\cite[Eq.~(1.7)]{KalvinCV}:
\begin{equation}\label{eq:C-calligraphic}
\begin{split}
 \mathcal C(\alpha)
 ={}&
 2\zeta_B'(0;\alpha,1,1)-2\zeta_R'(-1)
 -\frac{(\alpha-1)^2}{6\alpha}\log2\\
 &-\frac{\alpha-1}{12}+\frac12\log\alpha.
\end{split}
\end{equation}
Here $\zeta_B$ is the Barnes double zeta function and $\zeta_R$ is
the Riemann zeta function.

The general Aurell--Salomonson formula, rigorously proved in~\cite[Section~3.2]{KalvinJFA}, gives
\begin{equation}\label{eq:double-polygon-det}
\begin{split}
 \log\Det'\Delta_{\widehat P}
 ={}&
 \frac16\sum_{j=1}^N\sum_{\substack{i=1\\i\ne j}}^N
 \frac{(\alpha_i-1)(\alpha_j-1)}{\alpha_j}
 \log|p_i-p_j|
 -\sum_{j=1}^N\mathcal C(\alpha_j)
 \\
 &-\log C_{\boldsymbol\alpha,\boldsymbol p}^{\,2}
 -\zeta_{\boldsymbol\alpha}(0)
   \log\bigl(C_{\boldsymbol\alpha,\boldsymbol p}^{\,2}S\bigr)
 -\frac43\log2-4\zeta_R'(-1)
 +\frac16-\log\pi.
\end{split}
\end{equation}

For the Dirichlet Laplacian on $P$, the Aurell--Salomonson polygon
formula \cite[Eq.~(32)]{AurellSalomonsonFurther}, rewritten in the same
notation, is
\begin{equation}\label{eq:polygon-det}
\begin{split}
 \log\Det\Delta_{P,D}
 ={}&
 \frac1{12}\sum_{j=1}^N\sum_{\substack{i=1\\i\ne j}}^N
 \frac{(\alpha_i-1)(\alpha_j-1)}{\alpha_j}
 \log|p_i-p_j|
 -\frac12\sum_{j=1}^N\mathcal C(\alpha_j)
 \\
 &-\frac{\zeta_{\boldsymbol\alpha}(0)+1}{2}
  \log\bigl(C_{\boldsymbol\alpha,\boldsymbol p}^{\,2}S\bigr)
 +\frac14\sum_{j=1}^N\log\alpha_j
 \\
 &-2\zeta_R'(-1)+\frac1{12}
 -\frac16\log2-\frac12\log\pi.
\end{split}
\end{equation}
Unlike the closed-surface formula~\eqref{eq:double-polygon-det},
formula \eqref{eq:polygon-det} has not yet received a rigorous proof.
The derivation in \cite{AurellSalomonsonFurther} is partially
heuristic at the boundary corners.  Nevertheless, it is supported by
a number of explicitly evaluated examples and is generally believed
to be correct.  Here it is used as a model formula whose role in the
calculation is completely explicit.

Subtracting twice \eqref{eq:polygon-det} from
$\log(\Det'\Delta_{\widehat P}/S)$ cancels all terms involving the
prevertices $p_j$, the local functions $\mathcal C(\alpha_j)$, the
normalization $C_{\boldsymbol\alpha,\boldsymbol p}^{\,2}$, and the scale
$S$.  The remaining terms are
$-\log2-\frac12\sum_{j=1}^N\log\alpha_j$, so
\begin{equation}\label{eq:polygon-answer}
\boxed{\quad
 \cJ_{\rm BFK}(P)
 =
 \frac{\Det'\Delta_{\widehat P}/S}
      {(\Det\Delta_{P,D})^2}
 =
 \frac12\prod_{j=1}^N\alpha_j^{-1/2}.
\quad}
\end{equation}
This is precisely \eqref{eq:corner-conj}.  Although the two Laplacian
determinants separately depend on the moduli through the coordinates
$p_j$, their BFK quotient does not.  The corresponding jump determinant
is obtained by multiplying \eqref{eq:polygon-answer} by
$\length(\partial P)$, so it retains the perimeter and hence the scale.

For $N=3$, one can take $(p_1,p_2,p_3)=(-1,0,1)$.
Then \eqref{eq:double-polygon-det} and \eqref{eq:polygon-det} reduce
to the triangle formulas obtained from
\cite[Proposition~4.1]{KalvinCV} and \cite{AurellSalomonson},
respectively.

\section{A spindle cut into two spherical lunes}

The second example can be obtained directly from the spectrum.  The
determinant of the Laplacian on the rotationally
symmetric spindle was evaluated by Spreafico and Zerbini
\cite{SpreaficoZerbini} by separation of variables.  It was obtained
independently from the anomaly formula in
\cite[Section~3.1]{KalvinJFA}; the formula there applies to
constant-curvature spheres with two conical singularities and is not
restricted to antipodal singular points.

Let
\begin{equation}\label{eq:spindle-metric}
 g_{a,K}=K^{-1}\bigl(d\rho^2+a^2\sin^2\rho\,d\theta^2\bigr),
 \qquad
 0\leq\rho\leq\pi,\quad
 \theta\in\RR/2\pi\mathbb Z,
\end{equation}
where $a,K>0$.  This is a spherical spindle of curvature $K$ with cone
angle $2\pi a$ at each pole.  Its area is
\begin{equation}\label{eq:spindle-area}
 \Area(S^2_{a,K})=\frac{4\pi a}{K}.
\end{equation}

Cut the spindle along the two meridians $\theta=0$ and $\theta=\pi$.
The result consists of two congruent spherical lunes $L_{a,K}$.  Each
lune has two boundary corners of angle $\pi a$, and the total length of
the cutting curve is
\[
 \length(\Gamma)=\frac{2\pi}{\sqrt K}.
\]

Separation of variables gives the nonzero spectrum of the spindle:
\begin{equation}\label{eq:spindle-spectrum}
 \lambda_{m,n}
 =
 K\left(\frac{|m|}{a}+n\right)
  \left(\frac{|m|}{a}+n+1\right),
 \qquad
 m\in\mathbb Z,\quad n\in\mathbb N_0,
\end{equation}
with $(m,n)=(0,0)$ omitted.  The Dirichlet spectrum of one lune is
\begin{equation}\label{eq:lune-spectrum}
 \lambda^D_{m,n}
 =
 K\left(\frac{m}{a}+n\right)
  \left(\frac{m}{a}+n+1\right),
 \qquad
 m\in\mathbb N,\quad n\in\mathbb N_0.
\end{equation}
Indeed, the angular factors are $e^{im\theta}$ on the spindle and
$\sin(m\theta)$ on the lune.  Consequently their spectral zeta
functions satisfy the exact identity
\begin{equation}\label{eq:zeta-splitting}
 \zeta_{\rm sp}(s)
 =
 2\zeta_{L,D}(s)+Z_K(s),
 \qquad
 Z_K(s)=\sum_{n=1}^{\infty}[Kn(n+1)]^{-s}.
\end{equation}

We include the short evaluation of the last determinant.  For $K=1$,
binomial expansion followed by meromorphic continuation gives
\begin{equation}\label{eq:axis-zeta}
 Z_1(s)
 =
 \sum_{k=0}^{\infty}
 (-1)^k\frac{(s)_k}{k!}\zeta_R(2s+k).
\end{equation}
After meromorphic continuation, the $k=0$ and $k=1$ terms each give
$-1/2$ at $s=0$, so $Z_1(0)=-1$.  For the derivative,
\[
 Z_1'(0)
 =
 2\zeta_R'(0)-\gamma
 +\sum_{k=2}^{\infty}\frac{(-1)^k}{k}\zeta_R(k).
\]
The Taylor series of $\log\Gamma(1+z)$, evaluated at $z=1$ by an Abel
limit, gives
\[
 \sum_{k=2}^{\infty}\frac{(-1)^k}{k}\zeta_R(k)=\gamma.
\]
Therefore
\begin{equation}\label{eq:axis-values}
 Z_1(0)=-1,
 \qquad
 Z_1'(0)=-\log(2\pi).
\end{equation}
Since $Z_K(s)=K^{-s}Z_1(s)$,
\begin{equation}\label{eq:axis-det}
 Z_K'(0)=\log K-\log(2\pi),
 \qquad
 \exp[-Z_K'(0)]=\frac{2\pi}{K}.
\end{equation}

It now follows from \eqref{eq:zeta-splitting}--\eqref{eq:axis-det} that
\begin{equation}\label{eq:lune-answer}
\boxed{\quad
 \cJ_{\rm BFK}(L_{a,K})
 =
 \frac{\Det'\Delta_{S^2_{a,K}}/\Area(S^2_{a,K})}
      {(\Det\Delta_{L_{a,K},D})^2}
 =
 \frac{2\pi/K}{4\pi a/K}
 =
 \frac1{2a}.
\quad}
\end{equation}
Since the lune has two corners with $\alpha=a$,
\eqref{eq:lune-answer} is exactly the prediction of
Conjecture~\ref{conj:corner}.  Multiplication by
$\length(\Gamma)=2\pi/\sqrt K$ gives
$\Det_{\angle}'\cN_{\Gamma}=\pi/(a\sqrt K)$; thus the curvature cancels
only from the normalized quotient.

By contrast, cutting the same spindle along its smooth equator gives
the smooth factor $1/2$ \cite[Lemma~2.2]{KalvinCCM}.  The meridian cut
is the simplest one passing through conical points and hence carrying
corners.

\section{A spherical octant}

Let
\[
 O=\{(x,y,z)\in S^2:x,y,z\geq0\}
\]
be an octant of the unit round sphere.  It is a spherical triangle of
area $\pi/2$ with three angles equal to $\pi/2$.  Its mirror double
$\widehat O$ has area $\pi$ and three cone angles equal to $\pi$.
Equivalently,
\[
 \widehat O=S^2/H,
\]
where $H$ is the group of order four consisting of the even sign
changes of the three coordinates.

This is exactly the symmetric spherical double triangle appearing in
\cite[Appendix~A and Fig.~2]{KalvinCV}.  In the angle notation used
here,
\[
 \alpha_1=\alpha_2=\alpha_3=\frac12,
 \qquad
 |\boldsymbol\alpha|=\frac32,
 \qquad
 S=\Area(\widehat O)=\pi,
\]
and its curvature formula gives
\[
 K=\frac{2\pi(|\boldsymbol\alpha|-1)}{S}=1.
\]

The eigenvalues of the unit sphere are $l(l+1)$.  Dirichlet
eigenfunctions on $O$ are the spherical harmonics odd in each of
$x,y,z$.  Counting homogeneous harmonic polynomials of this parity
gives multiplicity $r$ in degree $l=2r+1$, $r\geq1$.  Hence
\begin{equation}\label{eq:octant-dirichlet-zeta}
 \zeta_{O,D}(s)
 =
 \sum_{r=1}^{\infty}
 r\bigl[(2r+1)(2r+2)\bigr]^{-s}.
\end{equation}
An $H$-invariant harmonic polynomial has all three coordinate parities
equal.  The odd sector is exactly \eqref{eq:octant-dirichlet-zeta},
whereas the even sector has multiplicity $r+1$ in degree $l=2r$.
After omitting the constant mode,
\begin{equation}\label{eq:octant-double-zeta}
 \zeta_{\widehat O}(s)
 =
 \sum_{r=1}^{\infty}
 (r+1)\bigl[2r(2r+1)\bigr]^{-s}
 +\zeta_{O,D}(s).
\end{equation}
Thus, with
\[
 \Phi(s)=\zeta_{\widehat O}(s)-2\zeta_{O,D}(s),
\]
we have
\begin{equation}\label{eq:octant-phi}
 \Phi(s)
 =
 \sum_{r=1}^{\infty}
 \left\{
 (r+1)\bigl[2r(2r+1)\bigr]^{-s}
 -
 r\bigl[(2r+1)(2r+2)\bigr]^{-s}
 \right\}.
\end{equation}

Binomial expansion, followed by the duplication identities for the
Hurwitz zeta and Barnes functions, gives
\begin{equation}\label{eq:octant-special-values}
\begin{aligned}
 \zeta_{O,D}(0)&=\frac{11}{48},
 &
 \zeta_{O,D}'(0)
 &=
 \frac12\zeta_R'(-1)
 +\frac12\log\pi+\frac14\log2-\frac1{16},
 \\
 \zeta_{\widehat O}(0)&=-\frac{13}{24},
 &
 \zeta_{\widehat O}'(0)
 &=\zeta_R'(-1)-\frac18.
\end{aligned}
\end{equation}
These values provide a direct check against \cite{KalvinCV}.  The
general formula there gives $-13/24$ for
$\alpha_1=\alpha_2=\alpha_3=1/2$, and Corollary~1.3 for these angles
and $S=\pi$ reduces to
\begin{equation}\label{eq:octant-cv-check}
 \log\Det'\Delta_{\widehat O}
 =
 \frac18-\zeta_R'(-1),
\end{equation}
which is exactly $-\zeta_{\widehat O}'(0)$ in
\eqref{eq:octant-special-values}.

It now follows from \eqref{eq:octant-special-values} that
\begin{equation}\label{eq:octant-phi-value}
 \Phi'(0)=-\log(\pi\sqrt2),
 \qquad
 \exp[-\Phi'(0)]=\pi\sqrt2.
\end{equation}
Since $\Area(\widehat O)=\pi$, it follows that
\begin{equation}\label{eq:octant-answer}
\boxed{\quad
 \cJ_{\rm BFK}(O)
 =
 \frac{\Det'\Delta_{\widehat O}/\Area(\widehat O)}
      {\bigl(\Det\Delta_{O,D}\bigr)^2}
 =
 \sqrt2.
\quad}
\end{equation}
This is precisely the prediction of Conjecture~\ref{conj:corner}:
\[
 \frac12\prod_{j=1}^3\left(\frac12\right)^{-1/2}
 =\sqrt2.
\]
The boundary of the octant consists of three great-circle arcs of
length $\pi/2$, so that $\length(\partial O)=3\pi/2$.  Therefore the
conjectural jump determinant itself is
\[
 \Det_{\angle}'\cN_{\partial O}
 =
 \length(\partial O)\,\cJ_{\rm BFK}(O)
 =
 \frac{3\pi}{\sqrt2}.
\]

For every spherical Coxeter triangle with angles
$(\pi/p,\pi/q,\pi/r)$, the same spectral calculation gives
$\cJ_{\rm BFK}=\frac12\sqrt{pqr}$, in agreement with
Conjecture~\ref{conj:corner}.  These calculations are lengthy but
otherwise analogous to the octant computation and are therefore
omitted.

\section{Consequence for constant-curvature polygons}

Let $P$ be a simply connected geodesic $N$-gon of constant curvature
$K\in\mathbb R$, whenever such a polygon exists, with interior angles
\[
 \pi\alpha_1,\ldots,\pi\alpha_N,
 \qquad 0<\alpha_j<2.
\]
Write $A=\Area(P)$ and $S=2A$.  The double $\widehat P$ is a
constant-curvature conical sphere of area $S$ with cone angles
$2\pi\alpha_j$.  Gauss--Bonnet gives
\begin{equation}\label{eq:constant-curvature-area}
 KA=\pi\left(\sum_{j=1}^N\alpha_j-(N-2)\right),
 \qquad
 KS=2\pi\left(\sum_{j=1}^N\alpha_j-(N-2)\right).
\end{equation}
Thus $\sum_j\alpha_j>N-2$ in the spherical case,
$\sum_j\alpha_j=N-2$ in the flat case, and
$\sum_j\alpha_j<N-2$ in the hyperbolic case.  When $K\ne0$ the area
is determined by the curvature and the angles; when $K=0$ the area is
an additional scale parameter.

Assuming Conjecture~\ref{conj:corner}, the determinant of the
Dirichlet Laplacian on $P$ is related to the determinant of the
Laplacian on its double by
\begin{equation}\label{eq:constant-curvature-polygon}
\boxed{\quad
 \Det\Delta_{P,D}
 =
 \left[
 2\left(\prod_{j=1}^N\alpha_j\right)^{1/2}
 \frac{\Det'\Delta_{\widehat P}}{S}
 \right]^{1/2}.
\quad}
\end{equation}

In the flat case this implication also runs in the reverse direction
from the model calculation in Section~3.  Once the corner BFK theorem
and \eqref{eq:corner-conj} have been proved independently, substituting
the rigorously established closed-double formula
\eqref{eq:double-polygon-det} into
\eqref{eq:constant-curvature-polygon} yields precisely
\eqref{eq:polygon-det}.  Thus the same theory would, in particular,
provide a rigorous proof of the Aurell--Salomonson formula for the
Dirichlet Laplacian on a Euclidean polygon.

The singular Polyakov anomaly formula of \cite{KalvinJFA} applies to
the closed double for arbitrary conical data and curvature.  It
compares determinants of conformally related conical metrics in terms
of the conformal factor but does not determine that factor; a closed
evaluation therefore requires an explicit uniformization.  In the flat
case it is supplied by the Schwarz--Christoffel data of Section~3,
whereas nonzero curvature requires solving the corresponding Liouville
equation.

For $N=3$ there are no moduli of the conical points.  The determinant
$\Det'\Delta_{\widehat P}$ is known explicitly as a function of the
angles and the area \cite{KalvinCV}: the required uniformization is
provided by the Schwarz triangle function.  Hence
\eqref{eq:constant-curvature-polygon} gives an explicit formula for the
Dirichlet determinant of every spherical, flat, or hyperbolic
geodesic triangle.  In this case, for $K\ne0$, equations
\eqref{eq:constant-curvature-area} and
\eqref{eq:constant-curvature-polygon} reduce extremal questions to
differentiating an explicit finite-dimensional function of the
angles.

For $N\ge4$, both Laplacian determinants in \eqref{eq:def-J} depend on
the moduli of the conical points, but this dependence cancels in their
BFK quotient.

For conical spheres obtained by gluing copies of a double triangle,
\cite{KalvinASNS} gives both an explicit Belyi uniformization and a
closed formula for the Laplacian determinant.  Consequently,
\eqref{eq:constant-curvature-polygon} is explicit whenever the polygon
double belongs to this Belyi-triangulated class.  In particular, the
determinants of the symmetric constant-curvature dihedra were evaluated
in \cite[Section~5.2]{KalvinASNS}.  Assuming
Conjecture~\ref{conj:corner}, the corresponding regular polygons are
stationary points of the Dirichlet determinant among
constant-curvature polygons with the same number of sides, area, and
curvature.

\section{Relation with jump, Loewner and Grunsky determinants}

It is useful to separate three related, but not identical, objects.

\paragraph{The spectral jump determinant.}
For a smooth curve on a closed Riemannian surface, the BFK quotient is
the zeta determinant of the Neumann jump operator divided by the length
of the curve.  Wang's determinant formula
\cite[Theorem~1.3]{Wang} expresses the Loewner energy of a smooth
Jordan loop \(\gamma\subset\mathbb{CP}^1\), relative to the reference
circle \(S^1\), as
\[
 I^L(\gamma)
 =
 12\left[
 \log\frac{\Det'\cN_\gamma}{\length(\gamma)}
 -
 \log\frac{\Det'\cN_{S^1}}{\length(S^1)}
 \right].
\]
The corresponding BFK quotient is
\begin{equation}\label{eq:gamma-bfk}
 \cJ_{\rm BFK}(\gamma,g)
 :=
 \frac{\Det'\Delta_{\mathbb{CP}^1,g}/\Area_g(\mathbb{CP}^1)}
      {\Det\Delta_{\Omega_+,D,g}\,
       \Det\Delta_{\Omega_-,D,g}}
 =
 \frac{\Det'\cN_\gamma}{\length(\gamma)}.
\end{equation}
For the reference circle, the conformal invariance of
\(\cJ_{\rm BFK}(S^1,g)\) and its value $1/2$ were obtained by the
present author in \cite[Lemma~2.10 and Section~3.4]{KalvinJFA},
respectively.  Consequently, Wang's determinant formula becomes
\begin{equation}\label{eq:wang-bfk}
 I^L(\gamma)
 =
 12\log\left(2\,\cJ_{\rm BFK}(\gamma,g)\right).
\end{equation}
Conjecture~\ref{conj:corner} therefore gives the following
BFK-renormalized polygonal analogue:
\begin{equation}\label{eq:wang-corner}
 12\log\left(2\,\cJ_{\rm BFK}(P)\right)
 =
 -6\sum_{j=1}^N\log\alpha_j.
\end{equation}
Here \(P\) is a geodesic polygon and the double is taken across its
boundary.  Since the ordinary Loewner energy is infinite at corners,
\eqref{eq:wang-corner} concerns the determinant side of Wang's formula;
comparison with reduced Loewner and Grunsky energies requires matching
the regularizations.

\paragraph{The planar interior--exterior problem.}
Wiegmann and Zabrodin \cite{WiegmannZabrodin} consider the sum of the
interior and exterior Dirichlet-to-Neumann operators of a planar
contour.  The determinant carries nontrivial conformal-welding
information and is related to the Fredholm
determinant involving a Neumann--Poincar\'e operator.  This is the natural object for the Dyson gas on a curved
contour.  It fits into the same general BFK mechanism after the plane
is compactified, but it should not be confused with the special
isometric mirror double.

\paragraph{Grunsky regularizations.}
For smooth Weil--Petersson curves, Loewner energy, the classical
Grunsky determinant and the jump determinant are tied together by
conformal anomaly and surgery formulas \cite{Wang,JohanssonViklund}.
For a curve with corners the Loewner energy is infinite.
Johansson--Viklund instead study a truncated Grunsky determinant and
identify the coefficient of its logarithmic divergence.  Their
conjectural reduced finite part is a natural candidate to compare with
Laplacian determinants, but such a comparison requires matching
regularizations. For the arc-Grunsky determinant in \cite{CourteautJohanssonViklund}, 
the corresponding Loewner-energy identity contains an explicit endpoint term. 
The authors also observe that the classical relation between the Grunsky and Neumann-Poincar\'e operators does not appear to extend to the arc-Grunsky setting.

The angle-dependent factor in \eqref{eq:corner-conj} is thus the local
factor in the conjectural BFK quotient.  The logarithmic divergence
found for truncated Grunsky matrices has a different form.  A complete
theory should explain exactly how the two regularizations are related.

\section{The conical-cut BFK problem}

The three model calculations in Sections~3--5 do not by themselves
construct $\Det_{\angle}'\cN$.  The missing result is a corner version
of the BFK theorem: a canonical determinant of the Neumann jump
operator for a cut allowed to pass through conical points, for which
the gluing identity \eqref{eq:smooth-bfk} remains valid.  In the
mirror-double case this determinant is expected to satisfy
\eqref{eq:corner-conj}.

In partiqular, such a theorem would turn the relations in Section~6 into determinant
formulas for Dirichlet Laplacians on constant-curvature polygons.  It should also clarify the relation with
the Loewner and Grunsky regularizations discussed in Section~7.
Polygonal cutting graphs, self-gluing, and Laplacians in bundles are
natural extensions of the same problem, but no general analytic
framework for them is proposed here.

\section*{Acknowledgements}

The author thanks Paul Wiegmann for discussions of a polygonal analogue
of the jump-determinant formula, beginning in 2023 and continuing more
recently in connection with the work of Johansson, Viklund and their
collaborators.  The author also thanks Werner M\"uller for a question
that prompted the present formulation.  Wiegmann and Zabrodin
acknowledge related earlier discussions with the author in
\cite{WiegmannZabrodin}.


\begin{thebibliography}{99}
\small
\setlength{\itemsep}{0pt}

\bibitem{AurellSalomonson}
E.~Aurell and P.~Salomonson,
\emph{On functional determinants of Laplacians in polygons and
simplicial complexes},
Comm. Math. Phys. \textbf{165} (1994), no.~2, 233--259;
\href{https://arxiv.org/abs/hep-th/9304031}{arXiv:hep-th/9304031}.

\bibitem{AurellSalomonsonFurther}
E.~Aurell and P.~Salomonson,
\emph{Further results on functional determinants of Laplacians in
simplicial complexes},
\href{https://arxiv.org/abs/hep-th/9405140}{arXiv:hep-th/9405140}.

\bibitem{BFK}
D.~Burghelea, L.~Friedlander and T.~Kappeler,
\emph{Mayer--Vietoris type formula for determinants of elliptic
differential operators},
J. Funct. Anal. \textbf{107} (1992), 34--65.

\bibitem{CourteautJohanssonViklund}
K.~Courteaut, K.~Johansson and F.~Viklund,
\emph{Planar Coulomb gas on a Jordan arc at any temperature},
\href{https://arxiv.org/abs/2504.19887}{arXiv:2504.19887}.

\bibitem{EdwardWu}
J.~Edward and S.~Wu,
\emph{Determinant of the Neumann operator on smooth Jordan curves},
Proc. Amer. Math. Soc. \textbf{111} (1991), no.~2, 357--363.

\bibitem{JohanssonViklund}
K.~Johansson and F.~Viklund,
\emph{Coulomb gas and the Grunsky operator on a Jordan domain with
corners},
Invent. Math. (2026),
\href{https://doi.org/10.1007/s00222-026-01417-5}
{doi:10.1007/s00222-026-01417-5};
\href{https://arxiv.org/abs/2309.00308}{arXiv:2309.00308}.

\bibitem{KalvinJFA}
V.~Kalvin,
\emph{Polyakov--Alvarez type comparison formulas for determinants of
Laplacians on Riemann surfaces with conical singularities},
J. Funct. Anal. \textbf{280} (2021), no.~7, Paper 108866;
\href{https://arxiv.org/abs/1910.00104}{arXiv:1910.00104}.

\bibitem{KalvinCCM}
V.~Kalvin,
\emph{Determinant of Friedrichs Dirichlet Laplacians on
$2$-dimensional hyperbolic cones},
Comm. Contemp. Math. \textbf{24} (2022), no.~10, Paper 2150107;
\href{https://arxiv.org/abs/2011.05407}{arXiv:2011.05407}.

\bibitem{KalvinCV}
V.~Kalvin,
\emph{Determinants of Laplacians for constant curvature metrics with
three conical singularities on $2$-sphere},
Calc. Var. Partial Differential Equations \textbf{62} (2023),
Paper 59;
\href{https://arxiv.org/abs/2112.02771}{arXiv:2112.02771}.

\bibitem{KalvinASNS}
V.~Kalvin,
\emph{Triangulations of singular constant curvature spheres via Belyi
functions and determinants of Laplacians},
preprint, 2023;
\href{https://arxiv.org/abs/2310.04882}{arXiv:2310.04882};
Ann. Sc. Norm. Super. Pisa Cl. Sci. \textbf{47} (2026),
published online 12 January 2026;
\href{https://journals.sns.it/index.php/annaliscienze/article/view/7204/2746}
{doi:10.2422/2036-2145.202501\_037}.

\bibitem{PeltolaWang}
E.~Peltola and Y.~Wang,
\emph{Large deviations of multichordal
$\mathrm{SLE}_{0+}$, real rational functions, and zeta-regularized
determinants of Laplacians},
J. Eur. Math. Soc. \textbf{26} (2024), no.~2, 469--535;
\href{https://arxiv.org/abs/2006.08574}{arXiv:2006.08574}.

\bibitem{SpreaficoZerbini}
M.~Spreafico and S.~Zerbini,
\emph{Spectral analysis and zeta determinant on the deformed spheres},
Comm. Math. Phys. \textbf{273} (2007), no.~3, 677--704;
\href{https://doi.org/10.1007/s00220-007-0229-z}
{doi:10.1007/s00220-007-0229-z};
\href{https://arxiv.org/abs/math-ph/0610046}
{arXiv:math-ph/0610046}.

\bibitem{Wang}
Y.~Wang,
\emph{Equivalent descriptions of the Loewner energy},
Invent. Math. \textbf{218} (2019), no.~2, 573--621;
\href{https://arxiv.org/abs/1802.01999}{arXiv:1802.01999}.

\bibitem{WiegmannZabrodin}
P.~Wiegmann and A.~Zabrodin,
\emph{Dyson gas on a curved contour},
J. Phys. A: Math. Theor. \textbf{55} (2022), no.~16, 165202;
\href{https://arxiv.org/abs/2111.09941}{arXiv:2111.09941}.

\end{thebibliography}
\end{document}